\documentclass[12pt]{article}

\usepackage{aeguill}
\usepackage[francais,english]{babel}
\usepackage[utf8]{inputenc}
\usepackage[T1]{fontenc}
\usepackage{amsfonts}
\usepackage{amssymb}
\usepackage{xspace}
\usepackage{graphicx}
\usepackage{amsmath}
\usepackage{amsthm}
\usepackage{stmaryrd}
\usepackage{MnSymbol}
\usepackage{endnotes}
\usepackage{yfonts}
\usepackage{makeidx}
\usepackage{multind}
\usepackage{fancyhdr}
\usepackage{pgf}
\usepackage{tikz}
\usepackage{dsfont}
\usepackage{bbm}
\usepackage[a4paper]{geometry}

\usetikzlibrary{positioning,arrows}

\theoremstyle{plain}
\newtheorem{thm}{Th\'eor\'eme}[section]
\newtheorem{prop}[thm]{Proposition}
\newtheorem{lem}[thm]{Lemme}

\theoremstyle{definition}
\newtheorem{df}[thm]{D\'efinition}

\newenvironment{dem}{\paragraph{Preuve}}{\qed\bigskip}

\newcommand{\N}{\mathbb{N}}
\newcommand{\Z}{\mathbb{Z}}

\newcommand{\R}{\mathbb{R}}

\newcommand{\Pb}{\mathbb{P}}

\newcommand{\sui}{_{n\geq 1}}
\newcommand{\fcar}{\mathbbm{1}}
\newcommand{\n[1]}{||#1||}
\newcommand{\tend}{\rightarrow}

\newcommand{\ab[1]}{|#1|}

\newcommand{\nti}{n\tend\infty}

\newcommand{\esp}{\mathbb{E}}

\newcommand{\Pj[1]}{\mathbb{P}(#1)}

\newcommand{\ma}{marche al\'eatoire\xspace}
\newcommand{\maX}{marche al\'eatoire sur $X$ associ\'ee \`a $G$ et $\mu$\xspace}

\DeclareMathOperator{\Supp}{Supp}

\begin{document}
\selectlanguage{francais}

\title{Un contre-exemple à la dichotomie récurrence/transience sur les espaces homogènes}
\author{Caroline Bruère}
\date{Juin 2016}
\maketitle

\selectlanguage{english}
\begin{abstract}
Take $G$ a locally compact second-countable group, 
and $H$ a subgroup of $G$. 
Choose $\mu$ a probability measure on $G$, 
such that the group spanned by its support is dense in $G$, 
and consider the Markov chain on the homogeneous space $X=G/H$
with transition probability $P_x=\mu *\delta_x$ for $x\in X$. 
Under some conditions on $G$, $H$, $\mu$, 
we know that this Markov chain is either everywhere recurrent 
or everywhere transient. 
A natural question is whether such a dichotomy is universally true.
The goal of this paper is to show it is not, 
even when $G$ is finitely generated, 
through the explicit construction of a counter example. 
The methods used include the study of stable laws, 
Gnedenko-Kolmogorov's local limit theorem for stable laws,
and the study of the time of first return to equilibrium. 
\end{abstract}

\selectlanguage{francais}

\begin{abstract}
Soit $G$ un groupe localement compact à base dénombrable,
$H$ un sous-groupe de $G$. 
Soit $\mu$ une mesure de probabilité sur $G$,
de support engendrant un sous-groupe dense dans $G$,
et considérons la chaîne de Markov sur l'espace homogène $X=G/H$
de probabilité de transition $P_x=\mu *\delta_x$ pour $x\in X$. 
Sous certaines conditions sur $G$, $H$, $\mu$,
on sait que cette chaîne de Markov est soit transiente en tout point, 
soit récurrente en tout point.
Une question naturelle est celle de savoir si une telle 
dichotomie est vraie universellement.
Le but de ce texte est de prouver que ce n'est pas le cas,
même quand $G$ est de type fini, 
à travers la construction explicite d'un contre-exemple. 
On utilisera l'étude du temps de premier retour à l'équilibre,
l'étude de lois stables,
et le théorème local-limite de Gnedenko-Kolmogorov 
pour les lois stables.  
\end{abstract}

\tableofcontents

\section{Introduction}

Considérons $G$ un groupe topologique localement compact à base dénombrable, 
$H$ un sous-groupe fermé de $G$, 
et $X=G/H$ le quotient de $G$ par $H$.
Soit $\mu$ une mesure de probabilité sur $G$.
La \emph{\maX} est la chaîne de Markov sur l'espace $X$ de probabilité de transition 
$P_x=\mu *\delta_x$ pour $x\in X$. 
Cette chaîne de Markov considérée est \emph{transitive} si
le support de $\mu$ engendre $G$ comme semi-groupe,
c'est à dire que pour tous points $x,\,y\in X$, 
la probabilité d'arriver en $y$ en temps fini en partant de $x$ est non-nulle).
On dit que $\mu$ est \emph{adaptée} si son support $\Supp \mu$ 
engendre un sous-groupe dense dans $G$. 
Notons $B=G^{\N^*}$, et $\beta=\mu^{\N^*}$ la mesure de probabilité produit sur $B$. 

\begin{df} \label{defrectrans}
On dit que la \ma sur $X$ est \emph{récurrente} en un point $x\in X$ 
s'il existe un compact $C$ de $X$ tel que 
\[\beta(\{b\in B\, |\, \forall n_0\in\N, \exists n\geq n_0 : b_n\cdots b_1x\in C\})=1.\] 
Elle est \emph{transiente} en un point $x\in X$ 
si pour tout compact $C$ de $X$, on a
\[\beta(\{b\in B\, |\, \exists n_0\in\N, \forall n\geq n_0 : b_n\cdots b_1x\notin C\})=1.\]
Elle est récurrente (respectivement transiente) 
sur tout $X$ si elle l'est en tout point. 
\end{df} 

On appelle \emph{théorème de dichotomie} 
une condition sur $G$, $H$, $\mu$ pour que la \maX associée
soit ou récurrente sur tout $X$, ou transiente sur tout $X$. 
Dans certains cas, ces théorèmes sont bien connus. 
Si $X$ est dénombrable, et que la \maX est transitive, alors elle est
soit récurrente sur tout $X$, soit transiente sur tout $X$ 
(voir le lemme \ref{equivrec}). 
Dans le cas où $X=G$, 
la dichotomie est toujours vraie
(voir par exemple \cite{GuiRaj11}).
Dans le cas où $H$ est un sous-groupe algébrique d'un groupe de Lie semi-simple, 
et $\mu$ est adaptée,
on a également un théorème de dichotomie, 
et une condition nécessaire et suffisante de récurrence (voir \cite{Bru1}). 
Le théorème de Hennion-Roynette (voir \cite{HennionRoynette})
donne un théorème de dichotomie 
quand la mesure $\mu$ est adaptée et étalée 
(par exemple si elle est absolument continue par rapport à la mesure de Haar),
et qu'il existe une mesure $\mu$-invariante sur $X$. 
Il est naturel de se demander si une telle dichotomie existe en général, 
pour $\mu$ adaptée,
sur les espaces homogènes de groupes localement compacts, à base dénombrable. 
Ce n'est pas le cas, même quand $G$ est de type fini. 

\begin{thm}\label{existeCE}
Il existe un groupe $G$ de type fini,
muni d'un sous-groupe $H$,
et de $\mu$ une mesure de probabilité sur $G$ adaptée 
tel que la \maX, en notant $X=G/H$,
ne soit ni transiente sur tout $X$, ni récurrente sur tout $X$. 
\end{thm}

Le but de ce texte est de prouver ce théorème,
en construisant un exemple explicite :
cette construction explicite fait l'objet de la partie \ref{contreexemple}. 
Nous commencerons par rappeler quelques propriétés des chaînes de Markov 
sur les espaces dénombrables, 
et des résultats sur les lois stables.
Nous étudierons ensuite des propriétés de la loi du premier temps de retour
à l'equilibre. 
Le lemme \ref{Zloistable}, qui caractérise la loi limite de la position
au premier temps de retour,
sera fondamental.
Dans la partie \ref{exZ2}, 
nous étudierons l'exemple classique de la marche simple sur $\Z^2$
à l'aide de ces résultats sur les temps de retour. 
Cette partie, qui n'est pas essentielle pour la suite du texte, 
ne consiste qu'en une exposition dans un cadre plus simple
des méthodes qui seront utilisées pour construire,
dans la partie \ref{contreexemple}, 
un contre-exemple à la dichotomie transience/récurrence.

\section{Une caractérisation de la récurrence}

Soit $Z$ un espace dénombrable discret. 
Soit $P$ un opérateur de Markov sur $Z$ ; 
il induit une chaîne de Markov sur $Z$,
dont l'espace des trajectoires est $Z^{\N}$. 
Pour $z\in Z$, notons $\Pb_z$ la mesure de probabilité sur l'ensemble des trajectoires issues de $z$
induite par $P$.
Notons $(X_n)\sui$ la suite des états de la chaîne de Markov induite par $P$. 
Pour un point $z\in Z$, pour $n\in\N$, notons $R_z^n$ le temps de $n$-me retour en $z$
de la chaîne de Markov : 
\[R^1_z=\inf\{k > 0\,|\,X_k=z\}\; ;\; R^n_z=\inf\{k>R^{n-1}_z\,|\,X_k=z\},\]
et $G_z$ le nombre de retour en $z$ :
\[G_z=\sum_{n=0}^{\infty}\fcar_{\{X_n=z\}}.\]
Rappelons les résultats suivants. 

\begin{lem}\label{equivrec} Soit $z\in Z$. Les assertions suivantes sont alors équivalentes :
\begin{enumerate}
\item Le point $z$ est récurrent ;
\item Il existe $y\in Z$ tel que $\forall n\in\N,\,\Pb_z(R^n_y<\infty)=1$ ;
\item Il existe $y\in Z$ tel que $\Pb_y(R^1_y<\infty)=1$ et $\Pb_z(R^1_y<\infty)=1$ ;
\item Il existe $y\in Z$ tel que $G_y=\infty$ $\Pb_z$-presque sûrement;
\item Il existe $y\in Z$ tel que $\esp_z(G_y)=\sum_{n=0}^{\infty} \Pb_z(X_n=y)=\infty$
et $\Pb_z(R^1_y<\infty)=1$.
\end{enumerate}

\end{lem}

\begin{dem} 
Par définition de l'espace $Z$, les compacts sont les ensembles finis de points. 
l'assertion $2$ est donc équivalent à la $1$ par définition de la récurrence.  
Soit $y\in Z$. 
En appliquant la propriété de Markov,
on obtient, pour tout $n\geq 1$, : 
\begin{align*}
\Pb_z(R^n_y<\infty)&=\Pb_z(R^1_y<\infty)\times\Pb_y(R^{n-1}_y<\infty)\\
&=\Pb_z(R^1_y<\infty)\times \Pb_y(R^1_y<\infty)^{n-1}. 
\end{align*} 
Ainsi, l'assertion $3$ est équivalente à la $2$. 
De la même manière, la propriété de Markov permet d'écrire, 
pour $y\in Z$ et $k \geq 1$ : 
\begin{align*}
\Pb_z(G_y\geq k+1) & = \Pb_z(R^1_y<\infty)\times\Pb_y(G_y\geq k)\\
&=  \Pb_z(R^1_y<\infty)\times\Pb_y(R^1_y<\infty)\times \Pb_y(G_y\geq k-1)\\
&= \Pb_z(R^1_y<\infty)\times\Pb_y(R^1_y<\infty)^k.
\end{align*} 
Ainsi, si l'assertion $3$ est vraie, on a, pour tout $k\in \ N$, 
\[\Pb_z(G_y\geq k+1)=1,\]
et ainsi  
\[\Pb_z(G_y=\infty)=1,\]
et l'assertion $4$ est donc vérifiée. 
L'assertion $5$ est une conséquence immédiate de la $4$. 
Supposons à présent vraie l'assertion $5$. 
On a, d'après ce qui précède, l'expression suivante de $\esp_z[G_y]$ :
\begin{align*}
\esp_z[G_y] &= \sum_{k=1}^{\infty} k \Pb_z(G_y=k)\\
&= \sum_{k=1}^{\infty} \Pb_z(G_y\geq k)\\
&= \Pb_z(R^1_y<\infty) \times \sum_{k=1}^{\infty} \Pb_y(R^1_y<\infty)^{k-1}.
\end{align*}
La quantité $\esp_z[G_y]$ est donc infinie si et seulement si 
la probabilité $\Pb_y(R^1_y<\infty)$ est égale à $1$,
et l'assertion $3$ est donc vérifiée. 
\end{dem}

\begin{prop}\label{espGinfdonnerec} Soit $z\in Z$ vérifiant 
\[\esp_z[G_z]=\infty.\]
Le point $z$ est récurrent. 
\end{prop}

\begin{dem}
En utilisant les mêmes arguments que dans la démonstration du lemme \ref{equivrec}, 
on obtient 
\[
\esp_z[G_z] = \sum_{k=1}^{\infty} \Pb_z(R^1_z<\infty)^k,
\]
et donc 
\[ \Pb_z(R^1_z<\infty)= 1.\]
Par le lemme \ref{equivrec}, $z$ est récurrent. 
\end{dem}

\section{Lois stables}

Rappelons quelques faits sur les lois stables.
Soit $m$ une mesure de probabilité sur $\R$, 
et $F$ la \emph{fonction de distribution de $m$} : 
pour $x\in\R$, on pose $F(x)=m(]-\infty,\, x])$. 

\begin{df} La loi $m$ est \emph{stable} 
si pour tous $a_1,\,a_2>0$, pour tous $b_1,\,b_2\in\R$,
il existe $a>0$ et $b\in\R$ tels qu'on ait l'égalité, pour tout $x\in\R$ :
\[F_1*F_2=F',\]
où $F_1(x)=F(a_1x+b_1)$, $F_2(x)=F(a_2x+b_2)$, $F'(x)=F(ax+b)$. 
En d'autre termes, pour toutes variables aléatoires $X_1,\, X_2$ indépendantes
de loi $\mu$, 
toutes constantes $a_1,\,a_2>0$,  $b_1\in\R$, 
il existe $a>0,\, b\in\R$ telles que $a_1X_1+a_2X_2+b_1$ et $aX_1+b$ 
soient de même loi. 
\end{df}

Le théorème suivant montre tout l'intérêt de la notion de loi stable. 
Il est démontré dans \cite{KolGne54},\cite{Brei}, \cite{IbLin}. 

\begin{thm} La loi $m$ est stable si et seulement s'il existe
une suite $(X_k)_{k\in\N^*}$ 
de variables aléatoires à valeur réelles, indépendantes, identiquement distribuées, 
et des suites $(B_n)_{n\in\N^*}$ et $(A_n)_{n\in\N^*}$ de réels 
respectivement strictement positifs et quelconques 
telles que la suite
\[(\frac{X_1+\hdots+X_n}{B_n}-A_n)_{n\in\N^*}\]
converge en loi vers $m$. 
\end{thm}

La proposition suivante donne à la fois la définition de \emph{l'exposant d'une loi stable},
et une condition de convergence vers une loi stable d'exposant déterminé. 

\begin{prop}\label{DAloistable}
Soient $\mu$ une mesure de probabilité non-dégénérée 
de fonction de distribution $F$, 
et $(X_k)_{k\in\N^*}$ une suite de variables aléatoires à valeur réelles, 
indépendantes, de loi $\mu$. 
Alors il existe des suites $(B_n)_{n\in\N^*}$ et $(A_n)_{n\in\N^*}$ de réels 
respectivement strictement positifs et quelconques 
telles que la suite
\[(\frac{X_1+\hdots+X_n}{B_n}-A_n)_{n\in\N^*}\]
converge vers une loi stable 
si et seulement s'il existe deux constantes $c_+,\,c_-\geq 0$, non toutes deux nulles,
et une constante $\alpha$, $0<\alpha\leq 2$,
telles que  les propriétés suivantes sont vérifiées :
\begin{enumerate}
\item  
\[\lim_{x\tend +\infty} \frac{F(-x)}{1-F(x)}=\frac{c_-}{c_+} ;\]
\item si $c_+>0$, alors pour tout $a>0$, on a 
\[\lim_{x\tend +\infty} \frac{1-F(ax)}{1-F(x)}=\frac{1}{a^{\alpha}} ;\]
\item si $c_->0$, alors pour tout $a>0$, on a 
\[\lim_{x\tend +\infty} \frac{F(-ax)}{F(-x)}=\frac{1}{a^{\alpha}}.\] 
\end{enumerate}
On dit alors que $\mu$ est \emph{dans le domaine d'attraction 
d'une loi stable d'exposant $\alpha$}.
La loi stable $m$ vers laquelle converge, en loi, la suite ci-dessus, 
est dite \emph{d'exposant $\alpha$}. 
\end{prop}

Il est également intéressant d'avoir des informations sur la forme des coefficients $B_n$ et $A_n$ : c'est l'objet de la proposition \ref{formeAnBn}. 

\begin{prop}\label{formeAnBn}
Conservons les notations et hypothèses de la proposition \ref{DAloistable}.
Si $F$ est de la forme 
\[F(x)=\begin{cases}
		\frac{c_1+h_1(x)}{x^{\alpha}}\;\text{si}\; x<0\\
		1-\frac{c_2+h_2(x)}{x^{\alpha}}\;\text{si}\; x>0
		\end{cases},\]
avec $c_1,\,c_2\in\R$, $0<\alpha<2$, et $h_1,\,h_2$ de limites nulles en $\infty$,
$\mu$ est dans le domaine d'attraction 
d'une loi stable d'exposant $\alpha$, et les réels $B_n$ sont de la forme $cn^{1/\alpha}$,
pour un $c\in\R$. 
Si la loi $\mu$ est symétrique, alors on peut choisir $A_n=0$ pour $n\in\N^*$. 
\end{prop}

Pour les preuves des propositions \ref{DAloistable}
et \ref{formeAnBn}, 
on pourra voir \cite[Thm 9.34]{Brei}, \cite[Thm 2.6.7]{IbLin}, 
ou \cite[\S 25,\, \S 35]{KolGne54}.

\paragraph{Remarque} Les lois stables d'exposant $2$ sont les lois normales. 

\medskip

Le théorème local-limite pour les lois stables est dû à Gnedenko et Kolmogorov 
(voir par exemple \cite{KolGne54}) ; 
l'énoncé qui suit est le théorème \cite[Thm. 4.2.1]{IbLin}. 

\begin{thm}\label{TLLGneKol} (Théorème local-limite de Gnedenko-Kolmogorov)
Soient $\mu$ une mesure de probabilité de fonction de distribution $F$, 
et $(X_k)_{k\in\N^*}$ une suite de variables aléatoires à valeur entières, 
indépendantes, de loi $\mu$. 
On suppose $\mu$ portée par un réseau $\{a+kh\,|\,k\in\Z\}$,
avec $a\in\Z$, et $h\in\N$ maximal. 
On considère la suite des variables 
\[Z_n=X_1+\hdots+X_n.\]
Supposons que $\mu$ est dans le domaine d'attraction d'une loi stable 
d'exposant $0<\alpha\leq 2$.
Choisissons des suites $(B_n)_{n\in\N^*}$ et $(A_n)_{n\in\N^*}$ de réels 
respectivement strictement positifs et quelconques
telles que la suite de variables aléatoires $\frac{Z_n}{B_n}-A_n$
converge en loi vers 
une loi stable d'exposant $\alpha$, de densité notée $g$.
Alors on a 
\[\lim_{\nti} \sup_{k\in\Z} \ab[\frac{B_n}{h}\Pb(Z_n=an+kh)-g(\frac{an+kh}{B_n}-A_n)]=0.\]
\end{thm}

Ce théorème fournira un point clé de l'étude de nos deux exemples. 

\section{Temps de retour à l'équilibre}

\subsection{Lois des temps de retour à l'équilibre de la marche simple}

Dans ce texte, on utilisera la proposition \ref{espGinfdonnerec} 
pour montrer la récurrence d'un point. 
L'étude des temps de retour sera fondamentale pour étudier 
l'espérance du noyau de Green correspondant. 
Considérons la chaîne de Markov induite par l'action de $\Z$ sur lui-même,
via une mesure de probabilité $\mu$ sur $\Z$ symétrique à support compact. 
Notons $R_n$ le temps de n-me retour en $0$ 
d'une trajectoire partant de $0$.

\begin{lem}\label{loiR1} Il existe une constante $\tau_{\mu}>0$ telle que 
l'on dispose de l'équivalent suivant pour la loi de $R_1$ :
\[\Pb(R_1=n)\underset{\substack{\infty \\ n\, \text{pair}}}{\sim}\frac{\tau_{\mu}}{n^{3/2}}.\]

\end{lem}

Si $n$ est impair, on a bien sûr
\[\Pb(R_1=n)=0.\]

\begin{dem} Ce résultat est un cas particulier du théorème $8$ de \cite{Kes63}. 
\end{dem}

\subsection{Position au temps de retour en dimension 2}

Considérons l'action de $\Z^2$ sur lui-même 
et la mesure de probabilité 
$\mu=\frac{1}{4}(\delta_{1,\,1}+\delta_{1,\,-1}+\delta_{-1,\,1}+\delta_{-1,\,-1})$
sur $\Z^2$. 
On considère l'opérateur de moyenne sur $\Z^2$ associé $P_{\mu}$. 
Cet opérateur induit une chaîne de Markov sur $\Z^2$ ; on notera $S_k,\,T_k$ 
les variables aléatoires représentant les coordonnées d'une trajectoire au temps $k$,
et $\Pb$ la mesure de probabilité induite par $P_{\mu}$ 
sur l'espace des trajectoires partant de $O=(0,\,0)$. 
Notons $R_n$ le temps de $n$-me retour en $0$ de la suite $(T_k)_k$. 
Remarquons que la variable $R_n$ est indépendante des variables $(S_k)_k$. 
Notons $U_n=S_{R_n}$, et $Z_n=S_{R_n}-S_{R_{n-1}}$, pour $n\in\N^*$,
avec $Z=Z_1=S_{R_1}$. 
Les variables aléatoires $Z_n$ sont indépendantes, identiquement distribuées,
de loi qu'on notera $\nu$. 
Notons $F$ la fonction de distribution de $\nu$ : 
$F:m\mapsto \Pb(Z\leq m)$,
et
$\tilde{F} : m\mapsto \Pb(Z\geq m)$.
Montrons que $\nu$ est dans le domaine d'attraction d'une loi stable
d'exposant $1$. 

\begin{lem}\label{Zloistable} La loi $\nu$ de la variable aléatoire $Z=S_{R_1}$ est dans le domaine d'attraction d'une loi stable d'exposant $1$. 
Plus précisément, $\frac{U_n}{n}$ converge en loi vers une loi stable d'exposant $1$,
de densité $g$. 
\end{lem}

\begin{dem}

D'après la proposition \ref{formeAnBn},
$\nu$ étant symétrique,
il suffit de montrer
qu'on a l'égalité, pour une constante $\sigma>0$,
\[\lim_{+\infty}m(\tilde{F}(m))=\sigma>0.\]
Écrivons explicitement $\tilde{F}(m)$ pour $m\in\N^*$ :
\[\tilde{F}(m)=\sum_{\substack{k=m \\ k\, \text{pair}}}^{\infty}\Pb(S_k\geq m)\Pb(R_1=k)
=\sum_{\substack{k=m \\ k\, \text{pair}}}^{\infty}
\sum_{\substack{l=m \\ l\, \text{pair}}}^{k}\Pb(S_k=l)\Pb(R_1=k).\]
Le lemme \ref{loiR1} donne un équivalent de $\Pb(R_1=k)$, pour $k$ pair :
\[\Pb(R_1=k)\underset{\substack{\infty \\ k\, \text{pair}}}{\sim}\frac{\tau_{\mu}}{k^{3/2}}.\]
Le théorème local-limite \ref{TLLGneKol} donne un équivalent de $\Pb(S_k = l)$,
puisque $\mu$ est dans le domaine d'attraction d'une loi stable
d'exposant $2$
par le théorème central-limite :
\[\lim_{k\tend\infty} \sup_{\substack{l\in\Z \\ l=k \, \text{mod}\, 2 }} 
\ab[\sqrt{k}\Pb(S_k=l)-d_{\mu} e^{-l^2/2k}]=0,\]
pour une certaine constante $d_{\mu}>0$. 
Considérons la somme 
\[D_m=m\sum_{\substack{k=m \\ k\, \text{pair}}}^{\infty}
\sum_{\substack{l=m \\ l=k \,\text{mod}\, 2}} ^{k}\frac{\tau_{\mu}}{k^{3/2}} 
\frac{d_{\mu}}{\sqrt{k}} e^{-l^2/2k}
=\tau_{\mu}d_{\mu}\sum_{\substack{l=m \\ l\, \text{pair}}}^{\infty}
\sum_{\substack{k=l \\ k\, \text{pair}}}^{\infty}
\frac{m}{k^2} e^{-l^2/2k}.\]
Fixons $m\geq 4$, et fixons $l\geq m$ pair. 
On peut alors écrire 
\[\sum_{\substack{k=l \\ k\, \text{pair}}}^{\infty}
\frac{m}{k^2} e^{-l^2/2k}=
\sum_{p=\frac{l}{2}}^{\infty}
\frac{m}{4p^2} e^{-l^2/4p},\]
puis obtenir, par une comparaison série-intégrale : 
\[
\frac{m}{l^2} (1-e^{-l/2}) \leq \sum_{p=\frac{l}{2}}^{\infty}
\frac{m}{4p^2} e^{-l^2/4p} \leq \frac{m}{l^2} (1-e^{-l^2/(2l-4)})
\]
On en déduit
les inégalités, pour $m\geq 2$ fixé :
\begin{align*} m(1-e^{-m/2})\sum_{l=m}^{\infty} \frac{1}{l^2}
&\leq \sum_{\substack{l=m \\ l\, \text{pair}}}^{\infty}
\sum_{\substack{k=l \\ k\, \text{pair}}}^{\infty}
\frac{m}{k^2} e^{-l^2/2k}
\leq m  \sum_{l=m}^{\infty} \frac{1}{l^2} \\
1-e^{-m/2}
&\leq \sum_{\substack{l=m \\ l\, \text{pair}}}^{\infty}
\sum_{\substack{k=l \\ k\, \text{pair}}}^{\infty}
\frac{m}{k^2} e^{-l^2/2k}\leq \frac{m}{m-1} \\
\tau_{\mu}d_{\mu}(1-e^{-m/2}) & \leq D_m
\leq \tau_{\mu}d_{\mu}\frac{m}{m-1}.
\end{align*}
On obtient 
\[\lim_{m\tend\infty} D_m = \tau_{\mu}d_{\mu}>0.\]
Considérons à présent  la différence 
\[\Delta_m=\ab[m\tilde{F}(m)-D_m].\]
On a 
\[
\Delta_m \leq m \sum_{\substack{l=m \\ l\, \text{pair}}}^{\infty}
\sum_{\substack{k=l \\ k\, \text{pair}}}^{\infty}
\ab[ \Pb(S_k=l)(\Pb(R_1=k)-\frac{\tau_{\mu}}{k^{3/2}}) 
+ \frac{\tau_{\mu}}{k^{3/2}}(\Pb(S_k=l)- \frac{d_{\mu}}{\sqrt{k}}e^{-l^2/2k})]. \\
\]
Soit $\epsilon>0$. 
Il existe $m_0\in\N$ tel que pour tout $m\geq m_0$, pour tous $k,\,l\geq m$ pairs, on ait 
\begin{align*}
\ab[\Pb(R_1=k)-\frac{\tau_{\mu}}{k^{3/2}}]&<\frac{\epsilon}{k^{3/2}} ;\\
\ab[\Pb(S_k=l)- \frac{d_{\mu}}{\sqrt{k}}e^{-l^2/2k}]&<\frac{\epsilon e^{-l^2/2k}}{\sqrt{k}} ;\\
\Pb(S_k=l)& \leq \frac{2d_{\mu}e^{-l^2/2k}}{\sqrt{k}}.
\end{align*}
On en déduit, pour $m\geq m_0$, une majoration de $\Delta_m$ :
\begin{align*}
\Delta_m & \leq m \sum_{\substack{l=m \\ l\, \text{pair}}}^{\infty}
\sum_{\substack{k=l \\ k\, \text{pair}}}^{\infty}
 \frac{2d_{\mu}e^{-l^2/2k}}{\sqrt{k}}\frac{\epsilon}{k^{3/2}} 
+ \frac{\tau_{\mu}}{k^{3/2}}\frac{\epsilon e^{-l^2/2k}}{\sqrt{k}} \\
& \leq c_{\mu}m\epsilon   \sum_{\substack{l=m \\ l\, \text{pair}}}^{\infty}
\sum_{\substack{k=l \\ k\, \text{pair}}}^{\infty}\frac{e^{-l^2/2k}}{k^2} \\
& \leq c_{\mu}m\epsilon  \sum_{l=m}^{\infty} \frac{1}{l^2} \\
& \leq 4c_{\mu}\epsilon
\end{align*}
pour une certaine constante $c_{\mu}>0$. 
On obtient donc 
\[\lim_{+\infty}m(\tilde{F}(m))=\tau_{\mu}d_{\mu}>0,\]
le résultat recherché. 

\end{dem}

\paragraph{Remarque 1} La loi stable limite est une loi de Cauchy 
(c'est à dire qu'elle est d'exposant $1$ et symétrique par rapport à un réel $x_0$) 
dont on peut calculer la densité $g$ : $g(s)=\frac{1}{\pi(s^2+1)}$ 
(cf. \cite[Thm 9.27]{Brei} par exemple pour une expression 
explicite de la forme des lois stables).

\paragraph{Remarque 2} La loi $\nu$ ci-dessus a été étudiée 
dans des cadres plus généraux : 
dans le cas où $\mu$ a des moments d'ordre $2+\delta$, le théoreme 1.1 
de \cite{Uch10} en donne un équivalent asymptotique précis,
par exemple, 
qu'on aurait pu utiliser dans la démonstration du lemme \ref{Zloistable}.
Signalons que la loi de la position au premier temps de retour dans un ensemble 
est étudiée dans \cite{Spitz} et \cite{Kes87} par exemple.

\section{Exemple : récurrence de la marche simple sur $\Z^2$}\label{exZ2}

Les résultats de cette partie ne sont pas nécessaires pour la suite du texte ; 
le but est uniquement d'exposer les méthodes 
qui seront utilisées plus loin sur un exemple bien connu. 
Considérons l'action de $\Z^2$ sur lui-même et la mesure de probabilité $\mu$ sur $\Z^2$ : 
\[\mu=\frac{1}{4}(\delta_{1,\,1}+\delta_{1,\,-1}+\delta_{-1,\,1}+\delta_{-1,\,-1}),\] 
On considère l'opérateur de moyenne sur $\Z^2$ associé $P_{\mu}$. 
Cet opérateur induit une chaîne de Markov sur $\Z^2$ ; on notera $S_k,\,T_k$ 
les variables aléatoires représentant les coordonnées d'une trajectoire au temps $k$,
et $\Pb$ la mesure de probabilité induite par $P_{\mu}$ 
sur l'espace des trajectoires partant de $O=(0,\,0)$. 
Le résultat suivant est bien connu :

\begin{prop}\label{recZ2} Le point $O$ est récurrent.

\end{prop}

Notons $R_n$ le temps de $n$-me retour en $0$ de la suite $(T_k)_k$. 
Remarquons que la variable $R_n$ est indépendante des variables $(S_k)_k$. 
D'après la proposition \ref{espGinfdonnerec}, il suffit de prouver que l'espérance
\[\esp(G)=\sum_{n=0}^{\infty}\Pb(S_{R_n}=0)\]
 de la fonction
\[G=\sum_{n=0}^{\infty} \fcar_{\{S_{R_n}=0\}}\]
est infinie.
La proposition \ref{recZ2} se déduit immédiatement du lemme suivant. 

\begin{lem}\label{TLL2d}
Il existe $a>0$, $n_0\in\N$, tels que pour tout $n\geq n_0$, on a
\[\Pb(S_{R_n}=0)\geq \frac{a}{n}.\]
\end{lem}

\begin{dem}
D'après le lemme \ref{Zloistable}, la loi de $S_{R_n}$ 
est dans le domaine d'attraction d'une loi stable d'exposant $1$. 
En appliquant le théorème local-limite \ref{TLLGneKol} à cette loi,
on obtient le lemme \ref{TLL2d}.
\end{dem}

\section{Un espace homogène ni transient ni récurrent}\label{contreexemple}

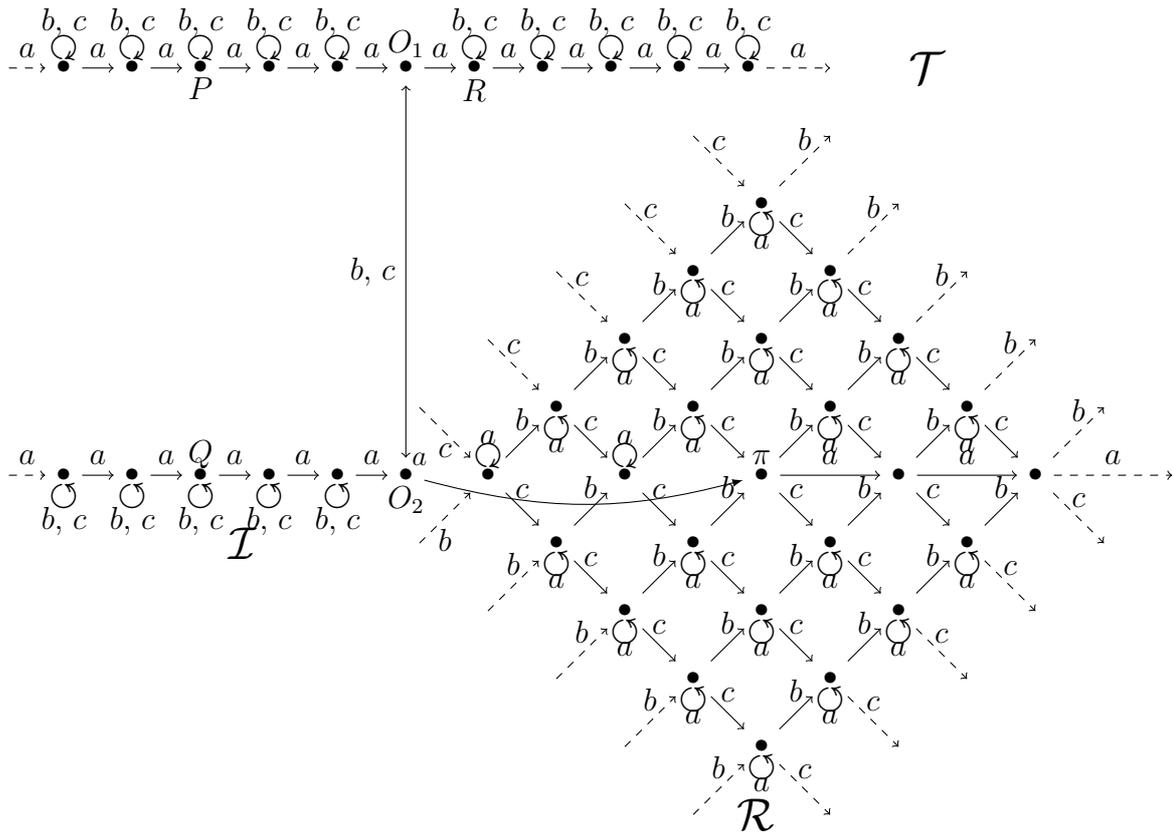
\begin{figure}

\centering

\begin{tikzpicture}[scale=0.9]

\node (h5a) at  (-5.2,6) {$\bullet$} ;
\node (h4a) at  (-4.2,6) {$\bullet$} ;
\node (h3a) at  (-3.2,6) {$\bullet$} ;
\node (h2a) at  (-2.2,6) {$\bullet$} ;
\node (h1a) at  (-1.2,6) {$\bullet$} ;
\node (h0) at  (-0.2,6) {$\bullet$} ;
\node (h1) at  (0.8,6) {$\bullet$} ;
\draw (-0.2, 6) node[above]{$O_1$} ;
\node (h2) at  (1.8,6) {$\bullet$} ;
\node (h3) at  (2.8,6) {$\bullet$} ;
\node (h4) at  (3.8,6) {$\bullet$} ;
\draw (0.8, 6) node[below]{$R$} ; 
\node (h5) at  (4.8,6) {$\bullet$} ;
\draw (-5.2, 6.3) node[above]{$b,\,c$} node{$\lcirclearrowdown$} ; 
\draw (-4.2, 6.3) node[above]{$b,\,c$} node{$\lcirclearrowdown$} ; 
\draw (-3.2, 6.3) node[above]{$b,\,c$} node{$\lcirclearrowdown$} ; 
\draw (-3.2, 5.7) node{$P$} ; 
\draw (-2.2, 6.3) node[above]{$b,\,c$} node{$\lcirclearrowdown$} ; 
\draw (-1.2, 6.3) node[above]{$b,\,c$} node{$\lcirclearrowdown$} ; 
\draw (0.8, 6.3) node[above]{$b,\,c$} node{$\lcirclearrowdown$} ; 
\draw (1.8, 6.3) node[above]{$b,\,c$} node{$\lcirclearrowdown$} ; 
\draw (2.8, 6.3) node[above]{$b,\,c$} node{$\lcirclearrowdown$} ; 
\draw (3.8, 6.3) node[above]{$b,\,c$} node{$\lcirclearrowdown$} ; 
\draw (4.8, 6.3) node[above]{$b,\,c$} node{$\lcirclearrowdown$} ; 
\draw [->](h5a) -- (h4a) node[midway,above]{$a$};
\draw [->](h4a) -- (h3a) node[midway,above]{$a$};
\draw [->](h3a) -- (h2a) node[midway,above]{$a$};
\draw [->](h2a) -- (h1a) node[midway,above]{$a$};
\draw [->](h1a) -- (h0) node[midway,above]{$a$};
\draw [->](h0) -- (h1) node[midway,above]{$a$} ;
\draw [->](h1) -- (h2) node[midway,above]{$a$};
\draw [->](h2) -- (h3) node[midway,above]{$a$};
\draw [->](h3) -- (h4) node[midway,above]{$a$};
\draw [->](h4) -- (h5) node[midway,above]{$a$};
\draw [->][dashed] (-6,6) -- (h5a) node[midway, above]{$a$} ;
\draw [->][dashed] (h5) -- (6,6) node[midway, above]{$a$}  ;
\node (b5a) at  (-5.2,0) {$\bullet$} ;
\node (b4a) at  (-4.2,0) {$\bullet$} ;
\node (b3a) at  (-3.2,0) {$\bullet$} ;
\node (b2a) at  (-2.2,0) {$\bullet$} ;
\node (b1a) at  (-1.2,0) {$\bullet$} ;
\node (b0) at  (-0.2,0) {$\bullet$} ;
\draw (-0.2, 0) node[below]{$O_2$} ;
\draw [->](b5a) -- (b4a) node[midway,above]{$a$};
\draw [->](b4a) -- (b3a) node[midway,above]{$a$};
\draw [->](b3a) -- (b2a) node[midway,above]{$a$};
\draw [->](b2a) -- (b1a) node[midway,above]{$a$};
\draw [->](b1a) -- (b0) node[midway,above]{$a$};
\draw (-5.2, -0.3) node[below]{$b,\,c$} node{$\circlearrowleft$} ; 
\draw (-4.2, -0.3) node[below]{$b,\,c$} node{$\circlearrowleft$} ; 
\draw (-3.2, -0.3) node[below]{$b,\,c$} node{$\circlearrowleft$} ; 
\draw (-3.2, 0.3) node{$Q$} ; 
\draw (-2.2, -0.3) node[below]{$b,\,c$} node{$\circlearrowleft$} ; 
\draw (-1.2, -0.3) node[below]{$b,\,c$} node{$\circlearrowleft$} ; 
\draw [->][dashed] (-6,0) -- (b5a) node[midway, above]{$a$} ;
\draw[<->] (b0) -- (h0) node[midway,left]{$b,\,c$} ;
\draw (7, 6) node{\Large$\;\;\;\;\;\mathcal{T}$} ; 
\draw (-3, -1) node{\Large$\;\;\;\;\;\mathcal{I}$} ; 
\draw (4.5, -5) node{\Large$\;\;\;\;\;\mathcal{R}$} ; 


\node (pi) at (5,0) {$\bullet$} ;
\node (x2y0) at (7,0) {$\bullet$} ;
\node (x4y0) at (9,0) {$\bullet$} ;
\node (x3ym1) at (8,-1) {$\bullet$} ;
\draw (8, -1.3) node[below]{$a$} node{$\circlearrowleft$} ;
\node (x1y1) at (6,1) {$\bullet$} ;
\draw (6, 0.7) node[below]{$a$} node{$\circlearrowleft$} ;
\node (x3y1) at (8,1) {$\bullet$} ;
\draw (8, 0.7) node[below]{$a$} node{$\circlearrowleft$} ;
\node (x2y2) at (7,2) {$\bullet$} ;
\draw (7, 1.7) node[below]{$a$} node{$\circlearrowleft$} ;
\node (xm1y1) at (4,1) {$\bullet$} ;
\draw (4,0.7) node[below]{$a$} node{$\circlearrowleft$} ;
\node (x0y2) at (5,2) {$\bullet$} ;
\draw (5,1.7) node[below]{$a$} node{$\circlearrowleft$} ;
\node (xm1y3) at (4,3) {$\bullet$} ;
\draw (4,2.7) node[below]{$a$} node{$\circlearrowleft$} ;
\node (x1y3) at (6,3) {$\bullet$} ;
\draw (6,2.7) node[below]{$a$} node{$\circlearrowleft$} ;
\node (xm2y0) at (3,0) {$\bullet$} ;
\draw (3, 0.3) node[above]{$a$} node{$\lcirclearrowdown$} ;
\node (xm4y0) at (1,0) {$\bullet$} ;
\draw (1, 0.3) node[above]{$a$} node{$\lcirclearrowdown$} ;
\node (xm3y1) at (2,1) {$\bullet$} ;
\draw (2,0.7) node[below]{$a$} node{$\circlearrowleft$} ;
\node (xm2y2) at (3,2) {$\bullet$} ;
\draw (3,1.7) node[below]{$a$} node{$\circlearrowleft$} ;
\node (x1ym1) at (6,-1) {$\bullet$} ;
\draw (6, -1.3) node[below]{$a$} node{$\circlearrowleft$} ;
\node (x2ym2) at (7,-2) {$\bullet$} ;
\draw (7, -2.3) node[below]{$a$} node{$\circlearrowleft$} ;
\node (x1ym3) at (6,-3) {$\bullet$} ;
\draw (6, -3.3) node[below]{$a$} node{$\circlearrowleft$} ;
\node (x0ym2) at (5,-2) {$\bullet$} ;
\draw (5, -2.3) node[below]{$a$} node{$\circlearrowleft$} ;
\node (xm1ym3) at (4,-3) {$\bullet$} ;
\draw (4, -3.3) node[below]{$a$} node{$\circlearrowleft$} ;
\node (xm1ym1) at (4,-1) {$\bullet$} ;
\draw (4, -1.3) node[below]{$a$} node{$\circlearrowleft$} ;
\node (xm3ym1) at (2,-1) {$\bullet$} ;
\draw (2, -1.3) node[below]{$a$} node{$\circlearrowleft$} ;
\node (xm2ym2) at (3,-2) {$\bullet$} ;
\draw (3, -2.3) node[below]{$a$} node{$\circlearrowleft$} ;
\node (x0y4) at (5,4) {$\bullet$} ;
\draw (5, 3.7) node[below]{$a$} node{$\circlearrowleft$} ;
\node (x0ym4) at (5,-4) {$\bullet$} ;
\draw (5, -4.3) node[below]{$a$} node{$\circlearrowleft$} ;
\draw (5,0) node[above]{$\pi$} ;

\draw[->, >=latex] (b0) to [bend right=15] (pi) node[very near end, above]{$a$} ;
\draw [->](pi) -- (x2y0) node[midway,above]{$a$};
\draw [->](x2y0) -- (x4y0) node[midway,above]{$a$};
\draw [->][dashed](x4y0) -- (11,0) node[midway,above]{$a$};

\draw[->][dashed](0, -1) -- (xm4y0) node[midway, below]{$b$} ; 
\draw[->](xm4y0) -- (xm3y1) node[midway, above]{$b$} ; 
\draw[->](xm3y1) -- (xm2y2) node[midway, above]{$b$} ; 
\draw[->](xm2y2) -- (xm1y3) node[midway, above]{$b$} ; 
\draw[->](xm1y3) -- (x0y4) node[midway, above]{$b$} ; 
\draw[->][dashed](x0y4) -- (6,5) node[midway, above]{$b$} ; 
\draw[->][dashed](1, -2) -- (xm3ym1) node[midway, above]{$b$} ; 
\draw[->](xm3ym1) -- (xm2y0) node[midway, above]{$b$} ; 
\draw[->](xm2y0) -- (xm1y1) node[midway, above]{$b$} ; 
\draw[->](xm1y1) -- (x0y2) node[midway, above]{$b$} ; 
\draw[->](x0y2) -- (x1y3) node[midway, above]{$b$} ; 
\draw[->][dashed](x1y3) -- (7,4) node[midway, above]{$b$} ; 
\draw[->][dashed](2,-3) -- (xm2ym2) node[midway, above]{$b$} ; 
\draw[->](xm2ym2) -- (xm1ym1) node[midway, above]{$b$} ; 
\draw[->](xm1ym1) -- (pi) node[midway, above]{$b$} ; 
\draw[->](pi) -- (x1y1) node[midway, above]{$b$} ; 
\draw[->](x1y1) -- (x2y2) node[midway, above]{$b$} ; 
\draw[->][dashed](x2y2) -- (8,3) node[midway, above]{$b$} ; 
\draw[->][dashed](3,-4) -- (xm1ym3) node[midway, above]{$b$} ; 
\draw[->](xm1ym3) -- (x0ym2) node[midway, above]{$b$} ; 
\draw[->](x0ym2) -- (x1ym1) node[midway, above]{$b$} ; 
\draw[->](x1ym1) -- (x2y0) node[midway, above]{$b$} ; 
\draw[->](x2y0) -- (x3y1) node[midway, above]{$b$} ; 
\draw[->][dashed](x3y1) -- (9,2) node[midway, above]{$b$} ; 
\draw[->][dashed](4,-5) -- (x0ym4) node[midway, above]{$b$} ; 
\draw[->](x0ym4) -- (x1ym3) node[midway, above]{$b$} ; 
\draw[->](x1ym3) -- (x2ym2) node[midway, above]{$b$} ; 
\draw[->](x2ym2) -- (x3ym1) node[midway, above]{$b$} ; 
\draw[->](x3ym1) -- (x4y0) node[midway, above]{$b$} ; 
\draw[->][dashed](x4y0) -- (10,1) node[midway, above]{$b$} ;


\draw[->][dashed](0, 1) -- (xm4y0) node[midway, below]{$c$} ; 
\draw[->](xm4y0) -- (xm3ym1) node[midway, above]{$c$} ; 
\draw[->](xm3ym1) -- (xm2ym2) node[midway, above]{$c$} ; 
\draw[->](xm2ym2) -- (xm1ym3) node[midway, above]{$c$} ; 
\draw[->](xm1ym3) -- (x0ym4) node[midway, above]{$c$} ; 
\draw[->][dashed](x0ym4) -- (6,-5) node[midway, above]{$c$} ; 
\draw[->][dashed](1, 2) -- (xm3y1) node[midway, above]{$c$} ; 
\draw[->](xm3y1) -- (xm2y0) node[midway, above]{$c$} ; 
\draw[->](xm2y0) -- (xm1ym1) node[midway, above]{$c$} ; 
\draw[->](xm1ym1) -- (x0ym2) node[midway, above]{$c$} ; 
\draw[->](x0ym2) -- (x1ym3) node[midway, above]{$c$} ; 
\draw[->][dashed](x1ym3) -- (7,-4) node[midway, above]{$c$} ; 
\draw[->][dashed](2,3) -- (xm2y2) node[midway, above]{$c$} ; 
\draw[->](xm2y2) -- (xm1y1) node[midway, above]{$c$} ; 
\draw[->](xm1y1) -- (pi) node[midway, above]{$c$} ; 
\draw[->](pi) -- (x1ym1) node[midway, above]{$c$} ; 
\draw[->](x1ym1) -- (x2ym2) node[midway, above]{$c$} ; 
\draw[->][dashed](x2ym2) -- (8,-3) node[midway, above]{$c$} ; 
\draw[->][dashed](3,4) -- (xm1y3) node[midway, above]{$c$} ; 
\draw[->](xm1y3) -- (x0y2) node[midway, above]{$c$} ; 
\draw[->](x0y2) -- (x1y1) node[midway, above]{$c$} ; 
\draw[->](x1y1) -- (x2y0) node[midway, above]{$c$} ; 
\draw[->](x2y0) -- (x3ym1) node[midway, above]{$c$} ; 
\draw[->][dashed](x3ym1) -- (9,-2) node[midway, above]{$c$} ; 
\draw[->][dashed](4,5) -- (x0y4) node[midway, above]{$c$} ; 
\draw[->](x0y4) -- (x1y3) node[midway, above]{$c$} ; 
\draw[->](x1y3) -- (x2y2) node[midway, above]{$c$} ; 
\draw[->](x2y2) -- (x3y1) node[midway, above]{$c$} ; 
\draw[->](x3y1) -- (x4y0) node[midway, above]{$c$} ; 
\draw[->][dashed](x4y0) -- (10,-1) node[midway, above]{$c$} ; 

\end{tikzpicture}
\caption{L'espace homogène $Z$ du groupe libre $F_3$}
\label{espaceZ}

\end{figure}

On va construire un espace homogène et un opérateur de Markov 
tels que la marche aléatoire associée n'est ni transiente, ni récurrente, 
prouvant ainsi le théorème \ref{existeCE}. 
Considérons le groupe libre à trois générateurs $G=F_3$ ; notons $a,\,b,\,c$ ces générateurs.
Considérons l'espace homogène $Z$ dessiné en figure \ref{espaceZ}.
L'espace $Z$ est un graphe comprenant trois parties :
\begin{itemize}
\item une partie $\mathcal{T}$ qui est une suite de points indexée par $\Z$,
\item une partie $\mathcal{I}$ qui est une suite de points indexée par $-\N$,
\item une partie $\mathcal{R}$ qui est un réseau indexé par 
$\{(i,\,j)\in\Z^2\,|\,i+j\,\text{pair}\,\}$.
\end{itemize}
On note $O_1$ un point de $\mathcal{T}$, $O_2$ l'extrémité de $\mathcal{I}$,
et $\pi$ un point de $\mathcal{R}$. 
L'action de $G$ est la suivante :
\begin{itemize}
\item les éléments $b$ et $c$ agissent par identité sur 
$\mathcal{T}\setminus \{O_1\}$ et sur $\mathcal{I}\setminus \{O_2\}$ ;
\item les éléments $b$ et $c$ échangent $O_1$ et $O_2$ ;
\item l'élément $a$ agit par translation sur $\mathcal{T}$ et $\mathcal{I}$ ;
\item l'élément $a$ envoie $O_2$ sur $\pi$ ;
\item en notant $\Delta$ le demi-axe horizontal de $\mathcal{R}$, d'extrémité $\pi$,
l'élément $a$ agit par identité sur les éléments de $\mathcal{R}\setminus\Delta$ ;
\item sur $\Delta$, l'élément $a$ agit par translation ;
\item l'élément $b$ agit par translation sur $\mathcal{R}$ 
le long de la première diagonale;
\item l'élément $c$ agit par translation sur $\mathcal{R}$ 
le long de la deuxième diagonale.
\end{itemize}
Le groupe $G$ agit transitivement sur $Z$ ; 
c'en est donc un espace homogène.

Considérons à présent la mesure de probabilité $\mu'$ à support compact suivante :
\[\mu'=\frac{1}{5}(\delta_a+\delta_b+\delta_{b^{-1}}+\delta_c+\delta_{c^{-1}}),\]
et la chaîne de Markov induite sur $Z$ par l'opérateur de moyenne $P_{\mu'}$. 
Etudions les trajectoires de cette chaîne de Markov.
Remarquons que l'action de $a$ est irréversible : 
graphiquement, quand on a parcouru un chemin dénoté $a$, 
on a une probabilité nulle de revenir en arrière. 
Ainsi, toutes 
les trajectoires issues du point $R$ (indiqué sur la figure \ref{espaceZ}) 
quittent tout compact avec probabilité $1$, par construction.
Considérons à présent les trajectoires issues du point $\pi$. 
Nous allons montrer qu'elles 
reviennent en $\pi$ avec probabilité $1$, et donc, d'après le lemme \ref{equivrec},
que le point $\pi$ est récurrent.

\begin{prop}\label{piestrec} Le point $\pi$ est récurrent.
\end{prop}

\begin{dem}

Notons $\Pb'$ la mesure de probabilité induite par $P_{\mu'}$ 
sur l'ensemble des trajectoires issues de $\pi$. 
Munissons $\mathcal{R}$ d'un système de coordonnées $(X,\,Y)$, 
représenté en figure \ref{espaceR}. 
Les points de $\mathcal{R}$ ont des coordonnées de la forme 
$(2k,\,2l)$ ou $(2k+1,\,2l+1)$, avec $k,\,l\in\Z$, 
et $\pi$ est de la forme $(0,\,0)$. 
Soit $(i,\,j)\in\mathcal{R}$. 
Les actions de $a,\,b,\,c$ s'écrivent alors :
\begin{itemize}
\item $b\cdot (i,\,j) = (i+1,\,j+1)$ ;
\item $c\cdot (i,\,j) = (i+1,\,j-1)$ ;
\item $a\cdot (0,\,j) = (0,\,j+2)$ lorsque $j\geq 0$;
\item $a\cdot (i,\,j) = (i,\,j)$ si $i\neq 0$ ou $j<0$. 
\end{itemize}
Notons $S'_k,\, T'_k$ les variables aléatoires représentant la position, 
respectivement horizontale et verticale, au temps $k$
d'une trajectoire issue de $\pi$, 
et $R'_n$ le temps de $n$-me retour en $0$ de $T'_k$.
D'après le lemme \ref{equivrec}, 
pour montrer que presque toutes les trajectoires issues de $\pi$ sont récurrentes, 
il suffit d'étudier la fonction $G'$,
\begin{equation}\label{formedeG}
G'=\sum_{n=0}^{\infty} \fcar_{\{S'_{R'_n}=0\}}.
\end{equation}
Remarquons que l'action du sous-groupe engendré par $b,\,c$ sur $\mathcal{R}$ 
est la même que celle décrite en partie \ref{exZ2}. 
Nous allons en fait étudier cette chaîne de Markov comme un décalage horizontal 
de la marche simple. 
Considérons l'action sur $\mathcal{R}$ du sous-groupe engendré par $b,\,c$,
muni de la mesure de probabilité $\mu$ équidistribuée :
\[\mu=\frac{1}{4}(\delta_b+\delta_{b^{-1}}+\delta_c+\delta_{c^{-1}}),\]
et la chaîne de Markov induite sur $\mathcal{R}$ par l'opérateur de moyenne $P_{\mu}$. 
Notons $\Pb$ la mesure de probabilité induite par $P_{\mu}$ 
sur l'ensemble des trajectoires issues de $\pi$,
et $S_k,\,T_k$ les variables aléatoires représentant la position, 
respectivement horizontale et verticale, au temps $k$ d'une trajectoire
issue de $\pi$. 
Notons enfin  $R_n$ le temps de $n$-me retour en $0$ de $T_k$.

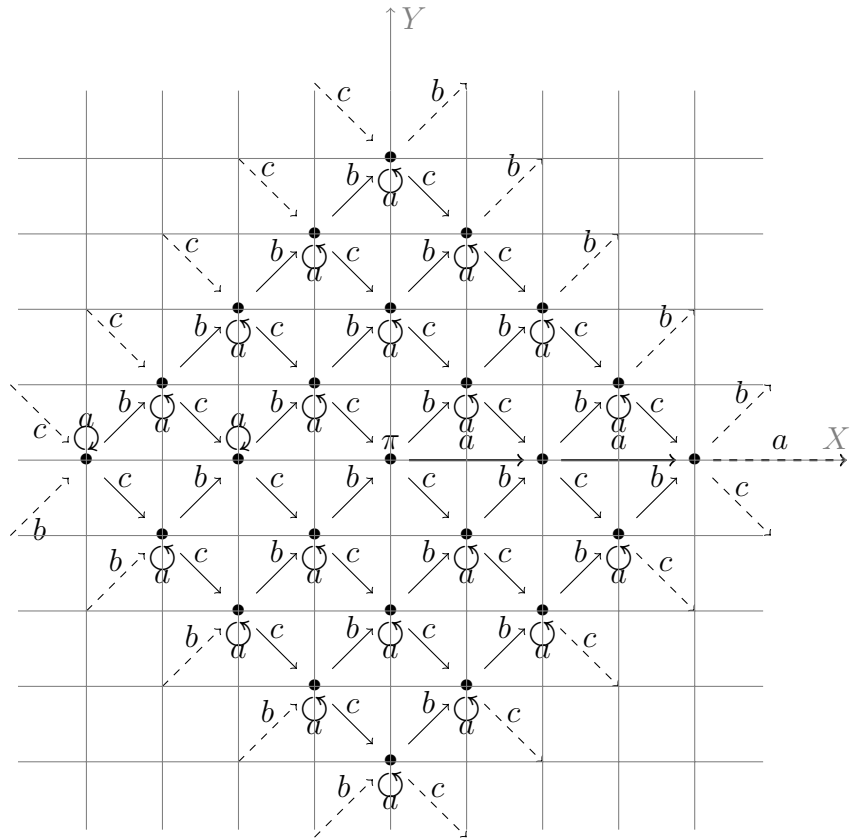
\begin{figure}
\centering
\begin{tikzpicture}

\node (pi) at (5,0) {$\bullet$} ;
\node (x2y0) at (7,0) {$\bullet$} ;
\node (x4y0) at (9,0) {$\bullet$} ;
\node (x3ym1) at (8,-1) {$\bullet$} ;
\draw (8, -1.3) node[below]{$a$} node{$\circlearrowleft$} ;
\node (x1y1) at (6,1) {$\bullet$} ;
\draw (6, 0.7) node[below]{$a$} node{$\circlearrowleft$} ;
\node (x3y1) at (8,1) {$\bullet$} ;
\draw (8, 0.7) node[below]{$a$} node{$\circlearrowleft$} ;
\node (x2y2) at (7,2) {$\bullet$} ;
\draw (7, 1.7) node[below]{$a$} node{$\circlearrowleft$} ;
\node (xm1y1) at (4,1) {$\bullet$} ;
\draw (4,0.7) node[below]{$a$} node{$\circlearrowleft$} ;
\node (x0y2) at (5,2) {$\bullet$} ;
\draw (5,1.7) node[below]{$a$} node{$\circlearrowleft$} ;
\node (xm1y3) at (4,3) {$\bullet$} ;
\draw (4,2.7) node[below]{$a$} node{$\circlearrowleft$} ;
\node (x1y3) at (6,3) {$\bullet$} ;
\draw (6,2.7) node[below]{$a$} node{$\circlearrowleft$} ;
\node (xm2y0) at (3,0) {$\bullet$} ;
\draw (3, 0.3) node[above]{$a$} node{$\lcirclearrowdown$} ;
\node (xm4y0) at (1,0) {$\bullet$} ;
\draw (1, 0.3) node[above]{$a$} node{$\lcirclearrowdown$} ;
\node (xm3y1) at (2,1) {$\bullet$} ;
\draw (2,0.7) node[below]{$a$} node{$\circlearrowleft$} ;
\node (xm2y2) at (3,2) {$\bullet$} ;
\draw (3,1.7) node[below]{$a$} node{$\circlearrowleft$} ;
\node (x1ym1) at (6,-1) {$\bullet$} ;
\draw (6, -1.3) node[below]{$a$} node{$\circlearrowleft$} ;
\node (x2ym2) at (7,-2) {$\bullet$} ;
\draw (7, -2.3) node[below]{$a$} node{$\circlearrowleft$} ;
\node (x1ym3) at (6,-3) {$\bullet$} ;
\draw (6, -3.3) node[below]{$a$} node{$\circlearrowleft$} ;
\node (x0ym2) at (5,-2) {$\bullet$} ;
\draw (5, -2.3) node[below]{$a$} node{$\circlearrowleft$} ;
\node (xm1ym3) at (4,-3) {$\bullet$} ;
\draw (4, -3.3) node[below]{$a$} node{$\circlearrowleft$} ;
\node (xm1ym1) at (4,-1) {$\bullet$} ;
\draw (4, -1.3) node[below]{$a$} node{$\circlearrowleft$} ;
\node (xm3ym1) at (2,-1) {$\bullet$} ;
\draw (2, -1.3) node[below]{$a$} node{$\circlearrowleft$} ;
\node (xm2ym2) at (3,-2) {$\bullet$} ;
\draw (3, -2.3) node[below]{$a$} node{$\circlearrowleft$} ;
\node (x0y4) at (5,4) {$\bullet$} ;
\draw (5, 3.7) node[below]{$a$} node{$\circlearrowleft$} ;
\node (x0ym4) at (5,-4) {$\bullet$} ;
\draw (5, -4.3) node[below]{$a$} node{$\circlearrowleft$} ;
\draw (5,0) node[above]{$\pi$} ;

\draw [->][thick](pi) -- (x2y0) node[midway,above]{$a$};
\draw [->][thick](x2y0) -- (x4y0) node[midway,above]{$a$};
\draw [->][thick][dashed](x4y0) -- (11,0) node[midway,above]{$a$};

\draw[->][dashed](0, -1) -- (xm4y0) node[midway, below]{$b$} ; 
\draw[->](xm4y0) -- (xm3y1) node[midway, above]{$b$} ; 
\draw[->](xm3y1) -- (xm2y2) node[midway, above]{$b$} ; 
\draw[->](xm2y2) -- (xm1y3) node[midway, above]{$b$} ; 
\draw[->](xm1y3) -- (x0y4) node[midway, above]{$b$} ; 
\draw[->][dashed](x0y4) -- (6,5) node[midway, above]{$b$} ; 
\draw[->][dashed](1, -2) -- (xm3ym1) node[midway, above]{$b$} ; 
\draw[->](xm3ym1) -- (xm2y0) node[midway, above]{$b$} ; 
\draw[->](xm2y0) -- (xm1y1) node[midway, above]{$b$} ; 
\draw[->](xm1y1) -- (x0y2) node[midway, above]{$b$} ; 
\draw[->](x0y2) -- (x1y3) node[midway, above]{$b$} ; 
\draw[->][dashed](x1y3) -- (7,4) node[midway, above]{$b$} ; 
\draw[->][dashed](2,-3) -- (xm2ym2) node[midway, above]{$b$} ; 
\draw[->](xm2ym2) -- (xm1ym1) node[midway, above]{$b$} ; 
\draw[->](xm1ym1) -- (pi) node[midway, above]{$b$} ; 
\draw[->](pi) -- (x1y1) node[midway, above]{$b$} ; 
\draw[->](x1y1) -- (x2y2) node[midway, above]{$b$} ; 
\draw[->][dashed](x2y2) -- (8,3) node[midway, above]{$b$} ; 
\draw[->][dashed](3,-4) -- (xm1ym3) node[midway, above]{$b$} ; 
\draw[->](xm1ym3) -- (x0ym2) node[midway, above]{$b$} ; 
\draw[->](x0ym2) -- (x1ym1) node[midway, above]{$b$} ; 
\draw[->](x1ym1) -- (x2y0) node[midway, above]{$b$} ; 
\draw[->](x2y0) -- (x3y1) node[midway, above]{$b$} ; 
\draw[->][dashed](x3y1) -- (9,2) node[midway, above]{$b$} ; 
\draw[->][dashed](4,-5) -- (x0ym4) node[midway, above]{$b$} ; 
\draw[->](x0ym4) -- (x1ym3) node[midway, above]{$b$} ; 
\draw[->](x1ym3) -- (x2ym2) node[midway, above]{$b$} ; 
\draw[->](x2ym2) -- (x3ym1) node[midway, above]{$b$} ; 
\draw[->](x3ym1) -- (x4y0) node[midway, above]{$b$} ; 
\draw[->][dashed](x4y0) -- (10,1) node[midway, above]{$b$} ;


\draw[->][dashed](0, 1) -- (xm4y0) node[midway, below]{$c$} ; 
\draw[->](xm4y0) -- (xm3ym1) node[midway, above]{$c$} ; 
\draw[->](xm3ym1) -- (xm2ym2) node[midway, above]{$c$} ; 
\draw[->](xm2ym2) -- (xm1ym3) node[midway, above]{$c$} ; 
\draw[->](xm1ym3) -- (x0ym4) node[midway, above]{$c$} ; 
\draw[->][dashed](x0ym4) -- (6,-5) node[midway, above]{$c$} ; 
\draw[->][dashed](1, 2) -- (xm3y1) node[midway, above]{$c$} ; 
\draw[->](xm3y1) -- (xm2y0) node[midway, above]{$c$} ; 
\draw[->](xm2y0) -- (xm1ym1) node[midway, above]{$c$} ; 
\draw[->](xm1ym1) -- (x0ym2) node[midway, above]{$c$} ; 
\draw[->](x0ym2) -- (x1ym3) node[midway, above]{$c$} ; 
\draw[->][dashed](x1ym3) -- (7,-4) node[midway, above]{$c$} ; 
\draw[->][dashed](2,3) -- (xm2y2) node[midway, above]{$c$} ; 
\draw[->](xm2y2) -- (xm1y1) node[midway, above]{$c$} ; 
\draw[->](xm1y1) -- (pi) node[midway, above]{$c$} ; 
\draw[->](pi) -- (x1ym1) node[midway, above]{$c$} ; 
\draw[->](x1ym1) -- (x2ym2) node[midway, above]{$c$} ; 
\draw[->][dashed](x2ym2) -- (8,-3) node[midway, above]{$c$} ; 
\draw[->][dashed](3,4) -- (xm1y3) node[midway, above]{$c$} ; 
\draw[->](xm1y3) -- (x0y2) node[midway, above]{$c$} ; 
\draw[->](x0y2) -- (x1y1) node[midway, above]{$c$} ; 
\draw[->](x1y1) -- (x2y0) node[midway, above]{$c$} ; 
\draw[->](x2y0) -- (x3ym1) node[midway, above]{$c$} ; 
\draw[->][dashed](x3ym1) -- (9,-2) node[midway, above]{$c$} ; 
\draw[->][dashed](4,5) -- (x0y4) node[midway, above]{$c$} ; 
\draw[->](x0y4) -- (x1y3) node[midway, above]{$c$} ; 
\draw[->](x1y3) -- (x2y2) node[midway, above]{$c$} ; 
\draw[->](x2y2) -- (x3y1) node[midway, above]{$c$} ; 
\draw[->](x3y1) -- (x4y0) node[midway, above]{$c$} ; 
\draw[->][dashed](x4y0) -- (10,-1) node[midway, above]{$c$} ;


\draw[thin, gray] (0.1, -4.9) grid (9.9, 4.9) ; 
\draw[->][gray] (5,4.9) -- (5,6) node[very near end, right]{$Y$} ;
\draw[->][gray] (9.9,0) -- (11,0) node[very near end, above]{$X$} ;

\end{tikzpicture}
\caption{Un système de coordonnées sur la partie $\mathcal{R}$}
\label{espaceR}
\end{figure}

Construisons à présent le décalage horizontal. 
Considérons la variable aléatoire $\eta$ de loi suivante : 
$\eta = 2m$ avec probabilité $\frac{4}{5^{m+1}}$, pour tout $m\in\N$. 
Notons $H_n=\sum_{k=0}^n \eta_k$ la somme de $n+1$ 
variables aléatoires indépendantes, identiquement distribuées, 
indépendantes des variables $S_k$ et $T_k$,
et de même loi que $\eta$.  
Remarquons que la loi de $\eta$ est à moments exponentiels.

D'après la proposition \ref{espGinfdonnerec}, 
pour montrer que $\pi$ est récurrent, il suffit de montrer que l'espérance 
par rapport à $\Pb'$ de
$G'=\sum_{n=0}^{\infty} \fcar_{\{S'_{R'_n}=0\}}$ 
est infinie. 
Or, par construction, on a pour $n\in\N$ :
\[\Pb'(S'_{R'_n}=0)=\Pb(S_{R_n}=-H_n)\]
On se ramène donc à étudier l'espérance par rapport à $\Pb$
de la fonction $G$ suivante :
\[G=\sum_{n=0}^{\infty} \fcar_{\{S_{R_n}=-H_n\}}.\]
On a :
\[\esp(G)=\sum_{n=0}^{\infty} \Pb(S_{R_n}=-H_n).\]
La proposition \ref{piestrec} se déduit alors immédiatement du lemme \ref{Gestinf},
via la proposition \ref{espGinfdonnerec}.
\end{dem}

\begin{lem}\label{Gestinf} On a l'égalité
\[\esp(G)=\infty.\]
\end{lem}

\begin{dem}
D'après le lemme \ref{Zloistable}, la loi de la variable $S_{R_1}$ est dans le domaine
d'attraction d'une loi stable d'exposant $1$, de densité $g$,
vers laquelle la suite $\frac{S_{R_n}}{n}$ converge en loi. 
Appliquons le théorème local-limite \ref{TLLGneKol} :
\[\lim_{\nti} \sup_{\substack{k\in\Z \\ k\, \text{pair}}} \ab[\frac{n}{2}
\Pb(S_{R_n}=k)-g(\frac{k}{n})]=0.\]
En notant $a=\sup_{t\in[-1,\,1]} g(t)>0$, 
il existe donc un $n_0\in\N^*$ tel que pour tout $n\geq n_0$, on ait
\[\forall k\in\llbracket -n,\,n \rrbracket,\, k\; \text{pair},\,\Pb(S_{R_n}=k)\geq \frac{a}{n}. \]
La variable aléatoire $\eta$ étant à moments exponentiels,
de moyenne $\frac{1}{2}$, 
on peut appliquer un principe des grandes déviations à $H_n =\sum_{k=1}^n \eta_k$ 
pour obtenir
l'existence $c>0$, $n_1\geq n_0$ tels que pour tout $n\geq n_1$, on ait
\[\Pb(H_n>n)\leq e^{-cn}.\]
Revenons au calcul de $\esp(G)$. 
Les variables $H_n$ et $S_{R_n}$ étant indépendantes, on peut écrire
\begin{align*}
\esp(G) & \geq \sum_{n=n_1}^{\infty} \sum_{k=0}^{n} \Pb(S_{R_n}=-k) \Pb(H_n=k) \\
& \geq \sum_{n=n_1}^{\infty}  \sum_{k=0}^{n} \frac{a}{n} \,  \Pb (H_n=k)\\
& \geq  \sum_{n=n_1}^{\infty}  \frac{a}{n} \, (1-e^{-cn})\\
& = \infty.
\end{align*}
\end{dem}

\subsection{Conclusion}

Ainsi, on a dans l'espace $Z$ :

\begin{itemize}
\item des points récurrents 
(par exemple le point $\pi$ ; c'est en fait le cas de tous les points de $\mathcal{R}$),
\item des points transients 
(par exemple le point $R$ indiqué sur la figure \ref{espaceZ}),
\item des points ni récurrents, ni transients (par exemple les points 
$P$, $Q$, $O_1$, $O_2$ indiqués sur la figure \ref{espaceZ}, 
ou tout point de $\mathcal{I}$).
\end{itemize}

L'espace $Z$, muni de l'action de $F_3$ 
et de la mesure de probabilité $\mu'$, 
n'est ni transient en tout point, ni récurrent en tout point :
le théorème \ref{existeCE} est donc vérifié. 
On n'a donc pas, en général, une dichotomie entre transience et récurrence 
sur les espaces homogènes. 

\nocite{Kes87}

\bibliography{recbib}{}

\begin{thebibliography}{10}

\bibitem{Brei}
L.~Breiman.
\newblock {\em Probability}.
\newblock Classics in Applied Mathematics. Society for industrial and applied
  mathematics, 1992.

\bibitem{Bru1}
C.~Bru{\`e}re.
\newblock Un crit{\`e}re de r{\'e}currence pour certains espaces homog{\`e}nes.
\newblock 2016.

\bibitem{GuiRaj11}
Y.~Guivarc'h and C.~Raja.
\newblock Polynomial growth, recurrence and ergodicity for random walks on
  locally compact groups and homogeneous spaces.
\newblock {\em Progress in Probability}, 64, 2011.

\bibitem{HennionRoynette}
H.~Hennion and B.~Roynette.
\newblock Un th\'eor\`eme de dichotomie pour une marche al\'eatoire sur un
  espace homog\`ene.
\newblock {\em Ast\'erisque}, 74:99--122, 1980.

\bibitem{IbLin}
I.~A. Ibragimov and Y.~V. Linnik.
\newblock {\em Independent and Stationary Sequences of Random Variables}.
\newblock Wolters-Noordhoff Publishing Groningen, 1971.

\bibitem{Kes63}
H.~Kesten.
\newblock Ratio theorems for random walks ii.
\newblock {\em Journal d'Analyse Math{\'e}matique}, 11(1):323--379, 1963.

\bibitem{Kes87}
H.~Kesten.
\newblock Hitting probabilities of random walks on $\mathbb{Z}^d$.
\newblock {\em Stochastic Processes and their Applications}, 1987.

\bibitem{KolGne54}
A.~N. Kolmogorov and B.~V. Gnedenko.
\newblock {\em Limit distributions for sums of independent random variables}.
\newblock Addison-Wesley, Cambridge, Mass., 1954.

\bibitem{Spitz}
F.~Spitzer.
\newblock {\em Principles of Random Walks}.
\newblock Springer.

\bibitem{Uch10}
K.~Uchiyama.
\newblock The hitting distributions of a line for two dimensional random walks
  2010.
\newblock {\em Transactions of the American Mathematical Society},
  362(5):2559--2588, 2010.

\end{thebibliography}
\bibliographystyle{abbrv}

\end{document}